\documentclass[11pt]{article}
\usepackage{amsmath}
\usepackage{amsthm}
\usepackage{amsfonts}

\newcommand{\mmp}{\mathbb{P}}

\newcommand{\od}{\overset{d}{=}}
\newcommand{\dod}{\overset{d}{\to}}
\newcommand{\tp}{\overset{P}{\to}}

\newcommand{\me}{\mathbb{E}}
\newcommand{\mr}{\mathbb{R}}
\newcommand{\mn}{\mathbb{N}}

\newcommand{\lin}{\underset{n\to\infty}{\lim}}

\newcommand{\lit}{\underset{t\to\infty}{\lim}}
\newenvironment{myproof}{\noindent {\it Proof} }{$\Box$ }
\newtheorem{thm}{Theorem}[section]
\newtheorem{lemma}[thm]{Lemma}

\newtheorem{cor}[thm]{Corollary}
\newtheorem{assertion}[thm]{Proposition}
\theoremstyle{definition}

\theoremstyle{remark}
\newtheorem{rem}[thm]{Remark}
\begin{document}
\title{On a random recursion related to absorption times of death Markov chains}
\date{October 30, 2007}
\author{Alex Iksanov\footnote{e-mail address:
iksan@unicyb.kiev.ua}\\ \small{\emph{Faculty of Cybernetics},
\emph{National T. Shevchenko University}},\\
\small{\emph{01033 Kiev, Ukraine}}\\
\ \\
Martin M\"ohle\footnote{e-mail address:
moehle@math.uni-duesseldorf.de}\\ \small{\emph{Mathematical
Institute}, \emph{University of D\"usseldorf}},\\
\small{\emph{40225 D\"usseldorf, Germany}}}

\maketitle

\begin{abstract}
   Let $X_1,X_2,\ldots$ be a sequence of random variables satisfying
   the distributional recursion $X_1=0$ and $X_n\od X_{n-I_n}+1$ for
   $n=2,3,\ldots$, where $I_n$ is a random variable with values in
   $\{1,\ldots,n-1\}$ which is independent of $X_2,\ldots,X_{n-1}$.
   The random variable $X_n$ can be interpreted as the absorption time
   of a suitable death Markov chain with state space $\mn:=\{1,2,\ldots\}$
   and absorbing state $1$, conditioned that the chain starts in the
   initial state $n$.

   This paper focuses on the asymptotics of $X_n$ as $n$ tends to
   infinity under the particular but important assumption that the
   distribution of $I_n$ satisfies $\mmp\{I_n=k\}=p_k/(p_1+\cdots+p_{n-1})$
   for some given probability distribution $p_k=\mmp\{\xi=k\}$, $k\in\mn$.

   Depending on the tail behaviour of the distribution of $\xi$, several
   scalings for $X_n$ and corresponding limiting distributions come into
   play, among them stable distributions and distributions of exponential
   integrals of subordinators.

   The methods used in this paper are mainly probabilistic. The key
   tool is a coupling technique which relates the distribution of
   $X_n$ to a random walk, which explains, for example, the appearance
   of the Mittag-Leffler distribution in this context.

   The results are applied to describe the asymptotics of the number of
   collisions for certain beta-coalescent processes.
\end{abstract}


\vspace{2mm}


\section{Introduction and main results} \label{intro}
Consider a death Markov chain $\{Z_k:k\in\mn_0:=\{0,1,\ldots\}\}$ with
state space $\mn:=\{1,2,\ldots\}$ and
transition probabilities $\pi_{ij}>0$ for $i,j\in\mn$ with $j<i$ and
$\pi_{ij}=0$ otherwise. For $n\in\mn$, define
$$
X_n\ :=\ \inf\{k\geq 1:Z_k=1\text{ given }Z_0=n\}.
$$
Note that $X_n\in\{1,2,\ldots,n-1\}$ almost surely.

Surprisingly, there seems to be very little known about the
asymptotic behavior of $X_n$ as $n$ tends to infinity. To our
knowledge, \cite{Cutsem} is one paper addressing this question.
However, the assumptions and the approach to be presented here are
completely different from those in \cite{Cutsem}.

The random variable $X_n$ can be interpreted as the
number of parts of the random composition $C_{n-1}$ of the integer
$n-1$, where the parts of the composition are (by definition) the
decrements of the Markov chain $\{Z_k: k\in \mn_0\}$. There are
several important articles in the
literature (\cite{BarGne, Gne, GnePitYor1, GnePitYor2}) with
asymptotic results on random compositions. However, in all these
papers the consistency of the random compositions for different
values of $n$ is a crucial assumption, i.e.~all these papers focus
on so called random composition structures. We do not assume this
consistency property here. Hence, our setting differs
significantly from that in the mentioned papers.

The key observation is that $X_n$ satisfies the distributional recursion
$X_1=0$ and
\begin{equation}\label{basic}
   X_n\ \od\ X_{n-I_n}+1,\quad n\in\{2,3,\ldots\},
\end{equation}
where $I_n$ is a random variable independent of $X_2,\ldots,X_{n-1}$
with distribution $\mmp\{I_n=k\}=\pi_{n,n-k}$, $k\in\{1,\ldots,n-1\}$.
The crucial assumption for the paper is that
\begin{equation}\label{basic2}
\mmp\{I_n=k\}\ =\ \frac{p_k}{p_1+\cdots+p_{n-1}},\qquad k,n\in\mn,
k<n,
\end{equation}
for some proper and non-degenerate probability distribution
\begin{equation}\label{distr}
   p_k\ :=\ \mmp\{\xi=k\},\quad k\in \mn, \ p_1>0.
\end{equation}
Throughout the paper $r(\cdot)\sim s(\cdot)$ means that
$r(\cdot)/s(\cdot)\to 1$ as the argument tends to infinity. The symbols
$\dod$, $\Rightarrow$, and $\tp$ denote convergence in law, weak
convergence, and convergence in probability, respectively, and
$X_n\dod (\Rightarrow, \tp)X$ means that the limiting relation holds
when $n\to\infty$. With $L$ we always denote a function slowly varying at
infinity.

Our main results given next are concerned with the limiting behaviour of
$X_n$ as $n\to\infty$. We begin with a weak law of large numbers.
\begin{thm}\label{weaklaw}
   If $\sum_{m=1}^n\sum_{k=m}^\infty p_k\sim L(n)$ for some function $L$
   slowly varying at $\infty$, then, as $n\to\infty$,
   \begin{equation}\label{weaklaw1}
       \dfrac{X_n}{\me X_n}\ \overset{P}{\to}\ 1
   \end{equation}
   and $\me X_n\sim n/L(n)$. In particular, if
   \begin{equation}\label{finitemean}
      m\ :=\ \me\xi\ <\ \infty,
   \end{equation}
   then $\me X_n\sim n/m$. If (\ref{finitemean}) holds, and if there
   exists a sequence of positive numbers $\{a_n:n\in\mn\}$ such
   that $X_n/a_n\tp 1$ as $n\to\infty$, then $a_n\sim n/m$.
\end{thm}
To formulate further results we need some more notation. For $C>0$
and $\alpha\in [1,2]$ let $\mu_\alpha$ be an $\alpha$-stable
distribution with characteristic function $\psi_\alpha(t)$, $t\in
\mr$ of the form
\begin{equation*}
\exp\{-|t|^\alpha
C\Gamma(1-\alpha)(\cos(\pi\alpha/2)+i\sin(\pi\alpha/2)\, {\rm
sgn}(t))\}, \ \ 1<\alpha<2;
\end{equation*}
\begin{equation*}
\exp\{-|t|C(\pi/2-i\log|t|\,{\rm sgn}(t))\}, \ \ \alpha=1;
\end{equation*}
\begin{equation*}
\exp(-(C/2)t^2), \ \ \alpha=2.
\end{equation*}
In the case when (\ref{finitemean}) holds, Theorem \ref{bas1}
provides necessary and sufficient conditions ensuring that $X_n$,
properly normalized and centered, possesses a weak limit.
\begin{thm}\label{bas1}
   If $m:=\me\xi<\infty$, then the following assertions are equivalent.
   \begin{enumerate}
      \item[(i)] There exist sequences of numbers $\{a_n,b_n:n\in\mn\}$
         with $a_n>0$ and $b_n\in\mr$ such that, as $n\to\infty$,
         $(X_n-b_n)/a_n$ converges weakly to a non-degenerate and
         proper probability law.
      \item[(ii)] Either $\sigma^2:=\mathbb{D}\xi<\infty$, or $\sigma^2=\infty$
         and for some $\alpha\in [1,2]$ and some function $L$ slowly
         varying at $\infty$,
         \begin{equation}\label{domain}
            \sum_{k=1}^n k^2p_k\ \sim\ n^{2-\alpha}L(n),\quad n\to\infty.
         \end{equation}
   \end{enumerate}
   If $\sigma^2<\infty$, then, with $b_n:=n/m$ and
   $a_n:=(m^{-3}C^{-1}\sigma^2n)^{1/2}$, the limiting law is $\mu_2$
   (normal with mean zero and variance $C$).\newline
   If $\sigma^2=\infty$ and (\ref{domain}) holds with $\alpha=2$, then, with
   $b_n:=n/m$ and $a_n:=m^{-3/2}c_n$, where $c_n$ is
   any sequence satisfying $\lim_{n\to\infty} nL(c_n)/c_n^2=C$, the limiting law
   is $\mu_2$.\newline
   If $\sigma^2=\infty$ and (\ref{domain}) holds with $\alpha\in [1,2)$, then,
   with $b_n:=n/m$ and $a_n:=m^{-(\alpha+1)/\alpha}c_n$,
   where $c_n$ is any sequence satisfying
   $$
   \lim_{n\to\infty}\dfrac{nL(c_n)}{c_n^\alpha}\ =\ \frac{\alpha}{2-\alpha}C,
   $$
   the limiting law is $\mu_\alpha$.
\end{thm}
\begin{rem}
%
   For $\sigma^2<\infty$, the same weak convergence result for $X_n$ was
   obtained in Theorem 4.1 in \cite{Cutsem} in a setting more general than
   ours. Note that 
   for $\alpha\in [1,2)$, (\ref{domain}) is equivalent to
   $\mmp\{\xi\geq n\}\sim (2-\alpha)n^{-\alpha}L(\alpha)/\alpha$,
   $n\to\infty$.
\end{rem}
If the mean of $\xi$ is infinite, the following Theorem
\ref{smallerone} (Theorem \ref{equalone}) points out conditions ensuring
that $X_n$, properly normalized (and centered), possesses a weak limit.
\begin{thm}\label{smallerone}
   Suppose that for some $\alpha\in(0,1)$ and some function $L$
   slowly varying at $\infty$
   \begin{equation}\label{regvar}
      \mmp\{\xi\geq n\}
      \ =\ \sum_{k=n}^\infty p_k\ \sim\ \dfrac{L(n)}{n^\alpha},\quad n\to\infty.
   \end{equation}
   Then, as $n\to\infty$,
   \begin{equation}\label{weak}
      \frac{L(n)}{n^\alpha}X_n\ \dod\ \int_0^\infty e^{-U_t}\,dt,
   \end{equation}
   where $\{U_t:t\geq 0\}$ is a subordinator with zero drift and L\'evy measure
   \begin{equation}\label{lev}
      \nu(dt)\ =\
   \frac{e^{-t/\alpha}}{(1-e^{-t/\alpha})^{\alpha+1}}\,dt,\quad t>0.
   \end{equation}
\end{thm}
It is instructive to present two proofs for Theorem \ref{smallerone},
namely a probabilistic proof and an analytic proof. The probabilistic
proof given in Section \ref{probab} reveals a relation between Eq.
(\ref{impo1}) and perpetuities. The analytic proof of Theorem
\ref{smallerone} presented in Section \ref{analytic} starts with the
distributional recursion (\ref{basic}), which implies that, for fixed
$k\in\mn$, the sequence $\{\me X_n^k: n\in \mn\}$ satisfies another
recursion. The structure of this last recursion permits a relatively
simple asymptotic analysis of $\me X_n^k$. In this way it is
possible to derive the convergence of the moments
$$
\lin
\me\left(\frac{L(n)}{n^\alpha}X_n\right)^k
\ =\
\me\left(\int_0^\infty e^{-U_t}\,dt\right)^k,\quad k\in\mn,
$$
which by a standard argument leads to (\ref{weak}).
\begin{thm}\label{equalone}
   Suppose that $\me\xi=\infty$ and that for some function $L$ slowly
   varying at $\infty$
   \begin{equation}\label{regvar1}
      \mmp\{\xi\geq n\}\ =\ \sum_{k=n}^\infty p_k\ \sim\ \dfrac{L(n)}{n}.
   \end{equation}
   Let $c$ be any positive function satisfying
   $\lim_{x\to\infty}xL(c(x))/c(x)=1$ and set $\psi(x):=x\int_0^{c(x)}
   \mmp\{\xi>y\}\,dy$. Let $b(x)$ be any positive function satisfying
   $$
   b(\psi(x))\ \sim\ \psi(b(x))\ \sim\ x,
   $$
   and set $a(x):=x^{-1}b(x)c(b(x))$. Then, $(X_n-b(n))/a(n)$
   converges weakly to the stable distribution $\mu_1$ with $C=1$.
\end{thm}
In the literature there exist two standard approaches to studying
distributional recursions. One approach is purely analytic and
based on a singularity analysis of generating functions (see,
for example, \cite{DIMR2,Pan2}). The other approach, called
{\it contraction method}, is more probabilistic (see
\cite{NeiRus,RosQu, RosRus}). It was remarked in \cite{IksMoe} that
recursions (\ref{basic}) which satisfy (\ref{basic2}) can be successfully
investigated by using probabilistic methods alone (completely different
from contraction methods). The present work extends ideas laid down in
\cite{IksMoe} for the particular case
$$
\mmp\{I_n=k\}\ =\ \dfrac{n}{n-1}\dfrac{1}{k(k+1)},
\quad k\in\{1,\ldots,n-1\}.
$$
The basic steps of the technique exploited can be summarized as follows.

Let $\xi_1,\xi_2,\ldots$ be independent copies of a random variable $\xi$
with distribution (\ref{distr}). Define $S_0:=0$, $S_n:=\xi_1+\cdots+\xi_n$
and $N_n:=\inf\{k\geq 1:S_k\geq n\}$, $n\in\mn$. Since $I_n\dod\xi$,
one may expect that the limiting behaviour of $X_n$ and $N_n$ is
similar, or at least that the limiting behaviour of the latter will
influence that of the former. To make this intuition precise, on the
probability space where $S_k$ and $N_n$ are defined,
we will construct (Section \ref{coupl}) random variables $M_n$
with the same distributions as $X_n$. Similarity in the limiting behaviour
of $M_n$ and $N_n$ is well indicated by asymptotic properties of
their difference. In particular, we will prove the
following.\newline (a) If $\me\xi<\infty$, then $M_n-N_n$ weakly
converges. Therefore, $M_n$, properly normalized and centered,
possesses a weak limit if and only if the same is true for $N_n$.\newline (b)
Assume now that $\me \xi=\infty$. (b1) If $\sum_{k=n}^\infty
p_k\sim L(n)/n$ and if $(N_n-b_n)/a_n$ weakly converges to
some $\mu$, then $(M_n-N_n)/a_n\tp 0$ which proves that
$(M_n-b_n)/a_n$ weakly converges to $\mu$. Thus in cases
(a) and (b1) a weak behaviour of $M_n$ and $N_n$ is the same. (b2)
If, for some $\alpha\in(0,1)$, $\sum_{k=n}^\infty p_k\sim
n^{-\alpha}L(n)$ and $N_n/a_n$ weakly converges to some
$\nu_1$, then $(M_n-N_n)/a_n$ weakly converges to some
$\nu_2$. Even though, the argument exploited above does not apply,
it will be proved that $M_n/a_n$ weakly converges to
$\nu_3\neq \nu_1$. Thus in this latter case a weak behaviour of
$M_n$ is not completely determined by that of $N_n$. Now
it is influenced by the weak behaviour of both $N_n$ and
$n-S_{N_n-1}$ to, approximately, the same extent. This observation
can be explained as follows. The probability of one big jump of
$S_n$ in comparison to cases (a) and (b1) is higher, and therefore
the epoch $N_n$ comes more "quickly". As a consequence, a
contribution to $M_n$ of the number of jumps in the sequence
$R_k^{(n)}$ (defined in Section \ref{coupl}), while $R_k^{(n)}$ is
travelling from $R_{N_n-1}^{(n)}=S_{N_n-1}$ to $n-1$, gets
significant.

It remains to review structural units of the paper not mentioned
so far. In Section \ref{nn} we investigate both the univariate and
the bivariate weak behaviour of $(N_n,n-S_{N_n-1})$, and discuss
their relation to exponential integrals of subordinators. Theorem
\ref{bas1}, \ref{weaklaw} and \ref{equalone} are proved in Section
4, 7 and 8 respectively. In Section 9 our main results apply to
the number of collisions in certain beta coalescent processes.
Possible generalizations of the results obtained and some directions
for future work are discussed in the final Section 10.
\section{A coupling}\label{coupl}
Fix $n\in\mn$. Define $R_0^{(n)}:=0$ and
$$
R_k^{(n)}\ :=\ R_{k-1}^{(n)}+\xi_k 1_{\{R_{k-1}^{(n)}+\xi_k<n\}},
\quad k\in\mn.
$$
Note that the sequence $\{R_k^{(n)}:k\in\mn_0\}$ is non-decreasing.
Let
$$
M_n\ :=\ \#\{i\in\mn:R_{i-1}^{(n)}\neq R_i^{(n)}\}
\ =\ \sum_{l=0}^\infty 1_{\{R_l^{(n)}+\xi_{l+1}<n\}}
$$
denote the number of jumps of the process $\{R_k^{(n)}:k\in\mn_0\}$.
Note that $M_1=0$ and that $1\leq M_n\leq n-1$ for $n\geq 2$. As
$p_1>0$, it follows from Lemma 1 in \cite{IksMoe} that the distribution
of $M_n$ satisfies the same recursion (\ref{basic}) as $X_n$. Hence,
the following lemma holds.
\begin{lemma}\label{xm}
   For each $n\in\mn$, the distribution of $M_n$ coincides with the
   distribution of the random variable $X_n$ introduced in Section \ref{intro}.
\end{lemma}
Fix $m,i\in\mn$. Define $\widehat{R}_0^{(m)}(i):=0$,
$$
\widehat{R}_k^{(m)}(i)\ :=\
\widehat{R}_{k-1}^{(m)}(i)+\xi_{i+k}1_{\{\widehat{R}_{k-1}^{(m)}(i)+\xi_{i+k}<m\}},
\quad k\in\mn,
$$
and
$$
\widehat{M}_n(i)\ :=\ \sum_{l=0}^\infty
1_{\{\widehat{R}_l^{(n)}(i)+\xi_{i+l+1}<n\}},
\quad n\in\mn_0.
$$
Our probabilistic proof of Theorem \ref{smallerone} relies upon the
following decomposition (\ref{impo1}).
\begin{lemma}\label{repr}
   For fixed $n\in\mn$ and any $i\in\mn$,
   \begin{equation}\label{eq1}
      \widehat{M}_n(i)\ \od\ M_n,
   \end{equation}
   and
   \begin{equation}\label{impo1}
      M_n-N_n+1
      \ =\ \widehat{M}_{n-S_{N_n-1}}(N_n)
      \ \od\ M^\prime_{n-S_{N_n-1}},
   \end{equation}
   where $\{M_n^\prime:n\in\mn\}$ has the same law as $\{M_n:n\in\mn\}$
   and is independent of $(N_n,n-S_{N_n-1})$.
\end{lemma}
\begin{proof}
   We have
   \begin{eqnarray*}
      M_n
      & = & \sum_{l=0}^\infty 1_{\{R_l^{(n)}+\xi_{l+1}<n\}}
      \ = \ \sum_{l=0}^{N_n-2}1+\sum_{l=N_n}^\infty 1_{\{R_l^{(n)}+\xi_{l+1}<n\}}\\
      & = & N_n-1+\sum_{l=0}^\infty
            1_{\{\widehat{R}_l^{(n-S_{N_n-1})}(N_n)+\xi_{N_n+l+1}<n-S_{N_n-1}\}}\\
      & = & N_n-1+\widehat{M}_{n-S_{N_n-1}}(N_n),
   \end{eqnarray*}
   and the first equality in (\ref{impo1}) follows. For any fixed $m\in\mn$,
   \begin{eqnarray*}
      &   & \hspace{-15mm}\mmp\{\widehat{M}_{n-S_{N_n-1}}(N_n)=m\}\\
      & = & \sum_{i=1}^n\sum_{j=0}^{n-1}
            \mmp\{\widehat{M}_{n-j}(i)=m, N_n=i,S_{N_n-1}=j\}\\
      & = & \sum_{i=1}^n\sum_{j=0}^{n-1}\mmp\{\sum_{l=0}^\infty
            1_{\{\widehat{R}_l^{(n-j)}(i)+\xi_{i+l+1}<n-j\}}=m, N_n=i,S_{N_n-1}=j\}.
   \end{eqnarray*}
   The sequence
   $\{\widehat{R}_l^{(n-j)}(i)+\xi_{i+l+1}: l\in\mn_0\}$ is
   independent of $1_{\{N_n=i, S_{N_n-1}=j\}}$ and has the same law
   as $\{(R_l^{(n-j)})^\prime+\xi^\prime_{l+1}: l\in \mn_0\}$, where
   $\{(R_l^{(\cdot)})^\prime: l\in \mn_0\}$ is constructed in the
   same way as the sequence without "prime" by using $\{\xi^\prime_k:
   k\in\mn\}$, an independent copy of $\{\xi_k: k\in\mn\}$. This
   implies (\ref{eq1}) and
   \begin{eqnarray*}
      &   & \hspace{-10mm}\mmp\{\widehat{M}_{n-S_{N_n-1}}(N_n)=m\}\\
      & = & \sum_{i=1}^n\sum_{j=0}^{n-1}
            \mmp\{\sum_{l=0}^\infty
            1_{\{(R_l^{(n-j)})^\prime+\xi^\prime_{l+1}<n-j\}}=m\}
            \mmp\{N_n=i,S_{N_n-1}=j\}\\
      & = & \mmp\{\sum_{l=0}^\infty
            1_{\{(R_l^{(n-S_{N_n-1})})^\prime+\xi^\prime_{l+1}<n-S_{N_n-1}\}}=m\}
      \ = \ \mmp\{M^\prime_{n-S_{N_n-1}}=m\},
   \end{eqnarray*}
   and the second equality in distribution in (\ref{impo1}) follows.
\end{proof}
\section{Results on $N_n$ and $n-S_{N_n-1}$: case $m=\infty$}\label{nn}
\subsection{Univariate results}
Below necessary and sufficient conditions are collected ensuring that
a properly normalized (without centering) $N_n$ weakly converges to
a non-de\-gener\-ate law (Proposition \ref{weakN}) and to $\delta_1$
(Proposition \ref{stab}).

We say that a random variable $\xi_\alpha$ has a Mittag-Leffler
distribution $\theta_\alpha$ with parameter $\alpha\in [0,1)$, if
$$
\me\xi_\alpha^n\ =\ \dfrac{n!}{\Gamma^n(1-\alpha)\,\Gamma(1+n\alpha)},
\quad n\in\mn.
$$
Note that the moments $\me\xi_\alpha^n$, $n\in\mn$, uniquely determine
the distribution. We also write $\theta_1$ for $\delta_1$.
\begin{assertion}\label{weakN}
   If (\ref{regvar}) holds for some $\alpha\in [0,1)$, then
   \begin{equation}\label{weakML}
      \dfrac{n^\alpha}{L(n)}N_n\ \Rightarrow\ \theta_\alpha.
   \end{equation}
   Conversely, assume that there exist positive real numbers $a(n)$,
   $n\in\mn$, such that $N_n/a(n)$ weakly converges to a
   non-degenerate law $\theta$. Then $a(n)\sim
   D\left(\sum_{k=n}^\infty p_k\right)^{-1}\sim Dn^\alpha/ L(n)$ for
   some constants $D>0$, $\alpha\in [0,1)$ and some function $L$ slowly
   varying at $\infty$, and (\ref{weakML}) holds.
\end{assertion}
\begin{rem}
   Proposition \ref{weakN} for $\alpha=0$ demonstrates that Theorem 6 in
   \cite{Feller} is wrong. For $\alpha\in (0,1)$, the implication (\ref{regvar})
   $\Rightarrow$ (\ref{weakML}) is well known (see, for example, Theorem 7 in
   \cite{Feller}). Our proof of Proposition \ref{weakN} seems to be new. It
   uses a technique introduced in \cite{DarKac} and simplified in \cite{BGT},
   Theorems 8.11.2 and 8.11.3. Note that $N_n$ is not the occupation time in the
   sense of Darling and Kac. Thus, before exploiting their approach, we
   had to prove Lemma \ref{DarKacle}, which is crucial for their technique
   to work.
\end{rem}
\begin{assertion}\label{stab}
   The following conditions are equivalent.
   \begin{enumerate}
      \item[(a)] $\sum_{m=1}^n\sum_{k=m}^\infty p_k\sim L(n)$ for some $L$ slowly
         varying at $\infty$.
      \item[(b)] $1-\sum_{n=1}^\infty e^{-sn}p_n\sim sL(1/s)$ as $s\downarrow 0$
         for some $L$ slowly varying at $\infty$.
      \item[(c)] The sequence $\{N_n: n \in \mn\}$ is relatively stable, i.e.
         there exist positive real numbers $a(n)$, $n\in\mn$, such that
         $N_n/a(n)\tp 1$.
   \end{enumerate}
   Moreover, if (a) holds, then $a(n)\sim\me N_n\sim n/L(n)$.
\end{assertion}
Put $P(s):=\sum_{n=1}^\infty e^{-sn}p_n$, $s\geq 0$, and
$h(s):=(1-P(s))^{-1}$, $s>0$. For $t\geq 0$ define $N_t:=\inf\{k\geq
1: S_k\geq t\}$. Then $N_t=N_1$ for $t\in [0,1]$, and $N_t=N_n$
for $t\in (n-1, n]$, $n=2,3,\ldots$
\begin{lemma}\label{DarKacle}
   Fix $k\in\mn$. Then, as $s\downarrow 0$,
   \begin{equation}\label{basLST}
      s\int_0^\infty e^{-st}\,\me N_t^k\,dt\ \sim\ k!\,h^k(s).
   \end{equation}
\end{lemma}
\begin{proof}
   For $k\in\{2,3,\ldots\}$ let $D_k$ denote the affine function of $k-2$
   positive variables of the form
   $$
   D_k(x_1,x_2,\ldots,x_{k-2})\ =\ \gamma_{0,k}+\sum_{i=1}^{k-2}\gamma_{i,k}x_i,
   $$
   with coefficients
   $\gamma_{i,k}\in\mr$, $i\in\{0,1,\ldots,k-2\}$. (These coefficients can
   be derived explicitly, but their exact values are of no use here.) For
   convenience, define $b_k(n):=\me N_n^k$, $k\in\mn$.
   We prove by induction on $k$ that
   \begin{equation}\label{recb}
      b_k(n)\ =\ c_k(n) + \sum_{i=1}^{n-1}b_k(n-i)\,p_i,\quad k\in\mn,
   \end{equation}
   with $c_1(n):=1$ and
   $$
   c_k(n)\ :=\ D_k(b_1(n),\ldots,b_{k-2}(n)) + k\,b_{k-1}(n),\quad k\geq 2.
   $$
   For $k=1$, Eq. (\ref{recb}) immediately follows from
   \begin{equation}\label{eqd}
      N_n\ \od\ 1+N^\prime_{n-\xi}1_{\{\xi<n\}},
      \quad n=2,3,\ldots,\quad N_1=1,
   \end{equation}
   where $\{N_n':n\in\mn\}$ is a copy of $\{N_n:n\in\mn\}$
   and $\xi$ is independent of $N_2',\ldots,N_{n-1}'$.
   Suppose (\ref{recb}) holds for $k\in\{1,2,\ldots,m-1\}$. Then,
   \begin{eqnarray*}
   &   & \hspace{-1cm}b_m(n)\ =\\
   & = & \sum_{i=0}^{m-2}
         \left(\hspace{-2mm}\begin{array}{c}m\\i\end{array}\hspace{-2mm}\right)
         \me(N_{n-\xi}'1_{\{\xi<n\}})^i
         +m\me(N_{n-\xi}'1_{\{\xi<n\}})^{m-1}
         +\me(N_{n-\xi}'1_{\{\xi<n\}})^m\\
   & = & 1 + m(b_1(n)-1) + \sum_{i=2}^{m-2}
         \left(\hspace{-2mm}\begin{array}{c}m\\i\end{array}\hspace{-2mm}\right)
         \big(b_i(n)-D_i(b_1(n),\ldots,b_{i-2}(n))\big)\\
   &   & \hspace{1cm}
         - mD_{m-1}(b_1(n),\ldots, b_{m-3}(n))
         + mb_{m-1}(n) + \sum_{i=1}^{n-1}b_m(n-i)p_i.
   \end{eqnarray*}
   The first four terms on the right-hand side form an affine
   function of $b_1(n)$, $\ldots$, $b_{m-2}(n)$, which implies
   (\ref{recb}) for $k=m$. Therefore, (\ref{recb}) is established.

   For $k\in\mn$ and $s>0$ define $B_k(s):=\sum_{n=1}^\infty
   e^{-sn}b_k(n)$ and $C_k(s):=\sum_{n=2}^\infty e^{-sn}c_k(n)$. Then,
   (\ref{recb}) is equivalent to
   \begin{equation}\label{impo}
      B_k(s)\ =\ \dfrac{e^{-s}+C_k(s)}{1-P(s)}
      \ =\ h(s)(e^{-s}+C_k(s)),\quad k\in\mn, s>0.
   \end{equation}
   We now verify by induction on $k$ that
   \begin{equation}\label{LST}
      s B_k(s)\ \sim\ k!\,h^k(s),\quad s\downarrow 0, k\in\mn.
   \end{equation}
   From $C_1(s)=e^{-2s}/(1-e^{-s})$, $s>0$, and
   (\ref{impo}) it follows that
   $$
   s B_1(s)
   \ =\ \dfrac{se^{-s}h(s)}{1-e^{-s}}
   \ \sim\ h(s),\quad s\downarrow 0.
   $$
   Thus, (\ref{LST}) holds for $k=1$. Suppose (\ref{LST}) holds for
   $k\in\{1,\ldots,m\}$ and check that
   $$
   s\,B_{m+1}(s)\ \sim\ (m+1)!\,h^{m+1}(s),\quad s\downarrow 0.
   $$
   The induction assumptions imply that
   $$
   s\sum_{i=2}^\infty e^{-si}D_{m+1}(b_1(i),\ldots,b_{m-1}(i))\ =\ o(h^m(s)),
   \quad s\downarrow 0.
   $$
   Therefore, by (\ref{impo}),
   \begin{eqnarray*}
      &   & \hspace{-3cm}sB_{m+1}(s)
      \ = \ sh(s)(e^{-s}+C_{m+1}(s))\\
      & = & sh(s)(e^{-s}+(m+1)(B_m(s)-e^{-s})\\
      &   & \hspace{1cm}
            +\sum_{i=2}^\infty e^{-si}D_{m+1}(b_1(i),\ldots,b_{m-1}(i)))\\
      & = & h(s)((m+1)sB_m(s)-mse^{-s}+o(h^m(s)))\\
      & \sim & (m+1)!\,h^{m+1}(s),\quad s\downarrow 0,
   \end{eqnarray*}
   and (\ref{LST}) is established. It remains to note that
   \begin{eqnarray*}
      &   & \hspace{-2cm}s\int_0^\infty e^{-st}\me N_t^k\,dt
      \ = \ s\sum_{j=1}^\infty \int_{j-1}^j e^{-st}\me N_j^k\,dt\\
      & = & (e^s-1)\sum_{j=1}^\infty e^{-sj}\me N_j^k
      \ \sim \ sB_k(s)
      \ \sim \ k!\,h^k(s),\quad s\downarrow 0.
   \end{eqnarray*}
\end{proof}
\begin{myproof}\emph{of Propositions \ref{weakN} and \ref{stab}}.
For $N$ sufficiently large, define $L(t):=L(n)$ for $t\in(n-1,n]$,
$n\in\{N,N+1,\ldots\}$. Then, (\ref{regvar}) is equivalent to
$\mmp\{\xi>x\}\sim x^{-\alpha}L(x)$, and condition (a) of Proposition
\ref{stab} is equivalent to $\int_0^x\mmp\{\xi>y\}\,dy\sim L(x)$.
Now, by Corollary 8.1.7 in \cite{BGT}, (\ref{regvar}) is equivalent to
$$
1-P(s)\ \sim\ \Gamma(1-\alpha)s^\alpha L(1/s),\quad s\downarrow 0,
$$
and conditions (a) and (b) of Proposition \ref{stab} are equivalent.
Regarding formally $\Gamma(0)$ as $1$, assume that $1-P(s)\sim
\Gamma(1-\alpha)s^\alpha L(1/s)$, $s\downarrow 0$, for some $\alpha
\in[0,1]$ or, equivalently,
$$
h(s)\ \sim\ \dfrac{1}{\Gamma(1-\alpha)s^\alpha L(1/s)},
\quad s\downarrow 0.
$$
We now proceed exactly as in the proof of Theorem 8.11.2 in
\cite{BGT}. Applying Karamata's theorem (\cite{BGT}, Theorem 1.7.6)
to (\ref{basLST}) gives
$$
\me N_t^k\ \sim\ \dfrac{k!}{\Gamma^k(1-\alpha)\Gamma(1+\alpha
k)}\dfrac{t^{\alpha k}}{L^k(t)}.
$$
Therefore,
\begin{equation}\label{mom1000}
   \lit\me\left(\dfrac{L(t)N_t}{t^\alpha}\right)^k
   \ =\ \dfrac{k!}{\Gamma^k(1-\alpha)\Gamma(1+\alpha k)},
   \quad k\in\mn,
\end{equation}
and, as $t\to\infty$, $L(t)N_t/t^\alpha\Rightarrow\theta_\alpha$,
which implies $L(n)N_n/n^\alpha\Rightarrow\theta_\alpha$.

Assume now that $N_n/a(n)\Rightarrow\theta$, and that either
$\theta=\delta_1$, or $\theta$ is non-degenerate. As the
sequence $\{N_n:n\in\mn\}$ is almost surely non-decreasing,
$\lim_{n\to\infty}N_n=\infty$ almost surely and $N_{n+1}\leq
N_n+1$ almost surely, we have $1\leq N_{n+1}/N_n\leq 1+1/N_n$
almost surely. Therefore, $\lim_{n\to\infty}N_{n+1}/N_n=1$ almost
surely and $\lim_{n\to\infty}a(n+1)/a(n)=1$. For $t>0$ define
$a(t):=a(n)$, $t\in(n-1,n]$. Then, by a sandwich argument,
$N_t/a(t)\Rightarrow\theta$ as $t\to\infty$.

If $\theta$ is non-degenerate, then from the proof of Theorem
8.11.3 in \cite{BGT} it follows that $a(t)\sim Dh(1/t)$ for some
$D>0$ and that the function $a$ regularly varies at $\infty$ with
exponent $\alpha\in[0,1)$. By Corollary 8.1.7,
$$
a(n)\ \sim\
\dfrac{D}{\Gamma(1-\alpha)\sum_{k=n}^\infty p_k}.
$$
Therefore, for some $\alpha\in[0,1)$, (\ref{regvar}) holds. By the
direct part of the proposition, (\ref{weakML}) holds as well.

If $\theta=\delta_1$, then we use a similar but simpler argument.
Let $T$ be an exponentially distributed random variable with mean
$1$ which is independent of $\{N_t:t\geq 0\}$. As in the proof of
Theorem 8.11.3 in \cite{BGT}, each sequence $r_n$ tending to $0$
contains a subsequence $\{s_n:n\in\mn\}$ satisfying $\lim_{n\to\infty}
s_n=0$, along which $\lim_{n\to\infty}a(t/s_n)/h(s_n)=f(t)$ at continuity
points of a non-decreasing function $f$. Therefore, $\lim_{n\to\infty}
a(T/s_n)/h(s_n)=f(T)$ almost surely. From (\ref{basLST}) it follows that
\begin{equation}\label{int15001}
   \lin\me\left(\dfrac{N_{T/s_n}}{h(s_n)}\right)^k\ =\ k!,
   \quad k\in\mn.
\end{equation}
Since $N_{T/s_n}/a(T/s_n)\tp 1$,
\begin{equation}\label{int1500}
   \dfrac{N_{T/s_n}}{h(s_n)}\ \tp\ f(T).
\end{equation}
Applying Fatou's lemma to (\ref{int15001}) with $k=1$ we conclude
that $f(T)<\infty$ almost surely. Also, (\ref{int15001}) implies
that, for each $k\in\mn$, the sequence $\{(N_{T/s_n}/h(s_n))^k:
n\in\mn\}$ is uniformly integrable which, in conjunction with
(\ref{int1500}), leads to $\me f^k(T)=k!$, $k\in\mn$. Since $\me
T^k=k!$, $k\in\mn$, and the sequence $\{k!:k\in\mn\}$ uniquely
determines (exponential) distribution, we conclude that $f(t)=t$,
$t>0$. The same argument as above can be repeated for any
sequence like $r_n$ which gives $a(t/s)/h(s)\to t$ as
$s\downarrow 0$ for each fixed $t>0$. Therefore, $a(t/s)/a(1/s)\to t$
as $s\downarrow 0$, which means that $a(t)\sim h(1/t)\sim t/L(t)$ as
$t\to\infty$ for some $L$ slowly varying at $\infty$. Hence,
$1-P(t)\sim tL(1/t)$ as $t\downarrow 0$.
\end{myproof}
\begin{rem}
Suppose (\ref{regvar}) holds for some $\alpha\in [0,1)$.
Then,
\begin{equation}\label{libeh}
   \lin\frac{L^k(n)}{n^{\alpha k}}\me N_n^k
   \ =\
   \dfrac{k!}{\Gamma^k(1-\alpha)\Gamma(1+\alpha k)},
   \quad k\in\mn.
\end{equation}
Suppose condition (a) of Proposition \ref{stab} holds. Then,
\begin{equation}\label{libeh1}
   \lin\frac{L^k(n)}{n^k}\me N_n^k=1,\quad k\in\mn.
\end{equation}
These observations immediately follow from (\ref{mom1000}). Note
also that (\ref{libeh}) is a particular case of Corollary 3.3
\cite{Port}.
\end{rem}
The next result is a corollary of Theorem \ref{weaklaw} and
Proposition \ref{stab}.
\begin{cor}\label{joint}
   Assume that (\ref{regvar1}) holds. Then, $\me N_n\sim\me M_n\sim
   n/m(n)$, where $m(x):=\int_0^x\mmp\{\xi>y\}\,dy$, $x>0$. Moreover,
   $$
   \dfrac{m(n)N_n}{n}\ \tp\ 1
   \quad\text{and}\quad
   \dfrac{m(n)M_n}{n}\ \tp\ 1.
   $$
   In particular, $M_n/N_n\tp 1$.
\end{cor}
\begin{proof}
   Condition (\ref{regvar1}) ensures that $m(x)$ belongs to
   the de Haan class $\Pi$, i.e.
   $\lim_{x\to\infty}(m(\lambda x)-m(x))/L(x)=\log\lambda$.
   In particular, $m(\cdot)$ is slowly varying at $\infty$. Since
   $\sum_{m=1}^n\sum_{k=m}^\infty p_k\sim m(n)$, Theorem
   \ref{weaklaw} together with Lemma \ref{xm} imply the
   result for $M_n$, and Proposition \ref{stab} implies the
   result for $N_n$.
\end{proof}
The next result is the key ingredient for our proof of Theorem
\ref{equalone}. Define $Y_n:=n-S_{N_n-1}$, $n\in\mn$.
\begin{assertion}\label{igrek}
   Assume that (\ref{regvar1}) holds. Then, for fixed $\delta>0$,
   \begin{equation}\label{rel}
      \me Y_n^\delta\ \sim\ \dfrac{n^\delta L(n)}{\delta\,m(n)},
\end{equation}
where $m(x):=\int_0^x\mmp\{\xi>y\}\,dy$, $x>0$. Furthermore, for
functions $a$ and $b$ as used in Theorem \ref{equalone},
\begin{equation}\label{rel2}
   \dfrac{b(n)Y_n}{n\,a(n)}\ \tp\ 0.
\end{equation}
\end{assertion}
\begin{proof}
   In the same way as in the proof of Proposition \ref{bivar}
   it follows that
   $$
   \me Y_n^\delta
   \ =\ \sum_{k=0}^{n-1}(n-k)^\delta\mmp\{\xi\geq n-k\}u_k,
   \quad n\in\mn,
   $$
   where $u_k:=\sum_{i=0}^k\mmp\{S_i=k\}$, $k\in\mn_0$. By
   Corollary \ref{joint}, $\me N_n\sim n/m(n)$. On the other
   hand, $\me N_n\sim\sum_{k=0}^n u_k$, $n\in\mn$. Thus,
   $\sum_{k=0}^n u_k\sim n/m(n)$ and, by Corollary 1.7.3 in
   \cite{BGT},
   $$
   U(s)
   \ :=\ \sum_{n=0}^\infty s^n u_n
   \ \sim\ \dfrac{1}{m((1-s)^{-1})(1-s)}
   \quad\text{as}\ s\uparrow 1.
   $$
   By the same Corollary
   $$
   V(s)
   \ :=\ \sum_{n=1}^\infty s^n n^\delta\mmp\{\xi\geq n\}
   \ \sim\ \dfrac{\Gamma(\delta)L((1-s)^{-1})}{(1-s)^\delta}
   \quad\text{as}\ s\uparrow 1.
   $$
   Therefore,
   $$
   \sum_{n=1}^\infty s^n\me Y_n^\delta
   \ =\ U(s)V(s)
   \ \sim\ \dfrac{\Gamma(\delta)}{(1-s)^{\delta+1}}
   \dfrac{L((1-s)^{-1})}{m((1-s)^{-1})}
   \quad\text{as}\ s\uparrow 1.
   $$
   The sequence $\{Y_n:n\in\mn\}$ is almost surely non-decreasing
   which implies that the sequence $\{\me Y_n^\delta: n\in\mn\}$
   is non-decreasing. Therefore, Corollary 1.7.3 in \cite{BGT} applies
   and proves (\ref{rel}).
   Recall that $\psi(x)=xm(c(x))$ and $c(x)\sim xL(c(x))$. Since
   $m(x)/L(x)\to\infty$, $c(x)\to\infty$ and
   $$
   \dfrac{\psi(x)}{c(x)}
   \ =\ \dfrac{xm(c(x))}{c(x)}
   \ \sim\ \dfrac{m(c(x))}{L(c(x))}
   $$
   as $x\to\infty$, we conclude that $\psi(x)/c(x)\to\infty$ as $x\to\infty$.
   Therefore,
   \begin{equation*}\label{j1}
      J_1(n)
      \ :=\ \dfrac{b([\psi(n)])}{a([\psi(n)])}
      \ =\ \dfrac{[\psi(n)]}{c(b([\psi(n)]))}
      \ \to\ \infty \quad \text{as}\ n\to\infty,
   \end{equation*}
   where $[x]$ denotes the integer part of $x$, and
   \begin{equation*}\label{j2}
      J_2(n)
      \ :=\ \dfrac{L([\psi(n)])}{m([\psi(n)])}\dfrac{b([\psi(n)])}{a([\psi(n)])}
      \ \sim\ \dfrac{L([nm(c(n))])}{m([nm(c(n))])}\dfrac{[nm(c(n))]}{nL(c(n))}
   \end{equation*}
   remains bounded for large $n$.

   Put $v(x):=xa(x)/b(x)=c(b(x))$. For fixed $\delta\in(0,1)$ and
   any $\epsilon>0$ we have, by Markov's inequality and by (\ref{rel}),
   $$
   \mmp\{Y_{[\psi(n)]}>v([\psi(n)]) \epsilon\}
   \leq \dfrac{\me Y^\delta_{[\psi(n)]}}{v^\delta([\psi(n)])\epsilon^\delta}
   \sim \dfrac{J_2(n)J^{\delta-1}_1(n)}{\delta\epsilon^\delta}
   \to 0 \
   \text{as} \ n\to\infty.
   $$
   The function $v$ is regularly varying at infinity with exponent $1$.
   Therefore, $\lim_{n\to\infty}v([\psi(n-1)])/v(\psi(n))=1$. Without loss
   of generality we can assume that $v$ is non-decreasing. If, for large $n$,
   $k\in([\psi(n-1)],[\psi(n)]]$, then
   $$
   \dfrac{Y_k}{v(k)}\ \leq\ \dfrac{Y_{[\psi(n)]}}{v([\psi(n-1)])}
   \quad\text{almost surely},
   $$
   and, by what we have already proved, as $n\to\infty$, the right-hand
   side tends to $0$ in probability, which proves (\ref{rel2}).
\end{proof}
\subsection{Some results on exponential integrals of
subordinators}\label{expsubord}
Let $\{Z_t:t\geq 0\}$ be a subordinator with zero drift which is
independent of $T$, an exponentially distributed random variable
with mean one. Set $Q:=\int_0^Te^{-Z_t}\,dt$, $M:=e^{-Z_T}$, and
$A:=\int_T^\infty e^{-Z_t}\,dt$. First of all, note that
\begin{eqnarray*}
   A_\infty
   & := & \int_0^\infty e^{-Z_s}\,ds
   \ = \ \int_T^\infty e^{-Z_s}\,ds\ +\ \int_0^T e^{-Z_s}\,ds\\
   & = & e^{-Z_T}\int_0^\infty e^{-(Z_{s+T}-Z_T)}\,ds
         \ +\ \int_0^T e^{-Z_s}\,ds.
\end{eqnarray*}
Therefore,
\begin{equation}\label{perp}
   A_\infty\ \od\ M A_\infty'+Q,
\end{equation}
where $A_\infty'$ is a copy of $A_\infty$ which is independent
of $(M,Q)$. The latter means that $A_\infty$ is a perpetuity
(see \cite{AlsIksRos} for the definition and recent results)
generated by the random vector $(M,Q)$. To verify (\ref{perp})
note that $\{Z_{s+t}-Z_t:s\geq 0\}$ is a subordinator which is
independent of $\{Z_v:v\leq t\}$ and has the same law as
$\{Z_u:u\geq 0\}$. Hence, for any Borel sets
$\mathcal{A}\in\mr^2$ and $\mathcal{B}\in\mr$,
\begin{eqnarray*}
   &   & \hspace{-2cm}\mmp\left\{
            (e^{-Z_T},A_T)\in \mathcal{A},\
            \int_0^\infty e^{-(Z_{s+T}-Z_T)}\,ds\in\mathcal{B}
         \right\}\\
   & = & \int_0^\infty e^{-t}\,
         \mmp\left\{
            (e^{-Z_t},A_t)\in\mathcal{A},\
            \int_0^\infty e^{-(Z_{s+t}-Z_t)}\,ds\in\mathcal{B}
         \right\}
         \,dt\\
   & = & \int_0^\infty e^{-t}\,
         \mmp\{(e^{-Z_t},A_t)\in\mathcal{A}\}\,dt\,
         \mmp\{A_\infty'\in\mathcal{B}\}\\
   & = & \mmp\{(e^{-Z_T}, A_T)\in\mathcal{A}\}\,
         \mmp\{A_\infty'\in\mathcal{B}\},
\end{eqnarray*}
and (\ref{perp}) follows.

Our next result generalizes Proposition 3.1 in \cite{Ca} dealing with
moments of $Q$, and a number of results concerning moments of
$\int_0^\infty e^{-Z_t}\,dt=Q+A$ (see, for example, Proposition 3.3
in \cite{Urbanik}).
\begin{assertion}\label{mixed}
   For $\lambda>0$ and $\mu\geq 0$
   $$
   \me Q^\lambda M^\mu
   \ =\ \dfrac{\lambda}{1+\varphi(\lambda+\mu)}\,\me Q^{\lambda-1}M^\mu,
   $$
   where $\varphi(s):=-\log\me e^{-sZ_1}$, $s\geq 0$. In particular,
   \begin{equation}\label{momentsQ}
      a_{n,m}\ :=\ \me Q^nM^m
      \ =\ \dfrac{n!}{\prod_{k=0}^n (1+\varphi(m+k))},
      \quad m,n\in \mn_0,
   \end{equation}
   $$
   b_{n,m}\ :=\ \me Q^nA^m
   \ =\ \dfrac{n!m!}{\prod_{k=0}^n (1+\varphi(m+k))\varphi(1)\cdots\varphi(m)},
   \quad m,n\in\mn_0.
   $$
   The moment sequences $\{a_{m,n}:m,n\in\mn_0\}$ and $\{b_{m,n}:m,n\in\mn_0\}$
   uniquely determine the laws of the random vectors $(M,Q)$ and $(A,Q)$
   respectively.
\end{assertion}
\begin{proof}
   For $t>0$ define $A_t:=\int_0^t e^{-Z_v}\,dv$. The following
   is essentially Eq. (3.1) in \cite{Ca}.
   $$
   A_t^\lambda e^{-\mu Z_t}\ =\ \lambda\int_0^t
   (A_t-A_v)^{\lambda-1}e^{-\mu(Z_t-Z_v)}e^{-(\mu+1)Z_v}\,dv.
   $$
   Since
   $$
   (A_t-A_v)^{\lambda-1}e^{-\mu(Z_t-Z_v)}
   \ =\ e^{-(\lambda-1)Z_v}
        \left(\int_0^{t-v}e^{-(Z_{s+v}-Z_v)}ds\right)^{\lambda-1}
        e^{-\mu(Z_t-Z_v)},
   $$
   and $\{Z_{s+v}-Z_v: s\geq 0\}$ is a subordinator which is independent of
   $\{Z_v:v\leq t\}$ and has the same law as $\{Z_t:t\geq 0\}$, we conclude
   that\newline
   $(\int_0^{t-v}e^{-(Z_{s+v}-Z_v)}\,ds)^{\lambda-1}
   e^{-\mu(Z_t-Z_v)}$ has the same law as $A_{t-v}^{\lambda-1}e^{-\mu
   Z_{t-v}}$ and is independent of $e^{-(\lambda-1)Z_v}$. Therefore,
   using Fubini's theorem,
   \begin{eqnarray*}
      \me A_T^\lambda e^{-\mu Z_T}
      & = & \int_0^\infty e^{-t}\me A_t^\lambda e^{-\mu Z_t}\,dt\\
      & = & \lambda\int_0^\infty e^{-t}
            \left(\int_0^t
               e^{-v\varphi(\lambda+\mu)}
               \me A_{t-v}^{\lambda-1}e^{-\mu Z_{t-v}}\,dv
            \right)\,dt\\
      & = & \lambda\int_0^\infty e^{-v\varphi(\lambda+\mu)}
            \left(\int_v^\infty
               e^{-t}\me A_{t-v}^{\lambda-1}e^{-\mu Z_{t-v}}\,dt
            \right)\,dv\\
      & = & \lambda\int_0^\infty e^{-v(\varphi(\lambda+\mu)+1)}\,dv
            \int_0^\infty e^{-u}\me A_u^{\lambda-1}e^{-\mu Z_u}\,du\\
      & = & \dfrac{\lambda}{1+\varphi(\lambda+\mu)}\,
            \me A_T^{\lambda-1}e^{-\mu Z_T}.
   \end{eqnarray*}
   Starting with
   \begin{equation}\label{momM}
      \me e^{-\mu Z_T}
      \ =\ \int_0^\infty e^{-t}\me e^{-\mu Z_t}\,dt
      \ =\ \int_0^\infty e^{-t(1+\varphi(\mu))}dt
      \ =\ \dfrac{1}{1+\varphi(\mu)},
   \end{equation}
   the formula for $a_{n,m}$ follows by induction. To prove that the
   law of $(M,Q)$ is uniquely determined by $\{a_{n,m}:n,m\in\mn_0\}$,
   it suffices to check that the marginal laws are uniquely determined
   by the corresponding moment sequences (see Theorem 3 in \cite{Peder}).
   Since $M\in[0,1]$ almost surely, the law of $M$ is trivially moment
   determinate. From (\ref{momentsQ}) it follows that
   $$
   \me Q^n\ =\ \dfrac{n!}{(1+\varphi(1)) \cdots (1+\varphi(n))},\quad n\in\mn.
   $$
   Set $f_n:=\me Q^n/n!$. The limit $f:=\lim_{n\to\infty} f_n/f_{n+1}$
   exists and is positive (it is finite, if $Z_t$ is compound
   Poisson, otherwise it is infinite). By the Cauchy-Hadamard
   formula, $f=\sup\{r>0:\me e^{rQ}<\infty\}$. Therefore, the law of
   $Q$ has finite exponential moments of some orders from which we deduce
   that this law is moment determinate.

   According to Proposition 3.3 in \cite{Urbanik}, $\me
   A_\infty^m=m!/(\varphi(1)\cdots\varphi(m))$, $m\in\mn_0$.
   In view of (\ref{perp}),
   \begin{eqnarray*}
      \me Q^nA^m
      & = & \me Q^nM^m \me A_\infty^m\\
      & = & \dfrac{n!\,m!}{\prod_{k=0}^n (1+\varphi(m+k))\varphi(1)\cdots\varphi(m)},
   \quad m,n\in\mn_0.
   \end{eqnarray*}
   In the same way as above for $(M,Q)$ it can be checked that the law
   of $(A,Q)$ is determined by the moment sequence. We omit the details.
\end{proof}
\subsection{A bivariate result}
Assume that (\ref{regvar}) holds, or, equivalently, that
\begin{equation}\label{regvar2}
   w(n)\ :=\ \dfrac{1}{\mmp\{\xi\geq n\}}
   \ =\ \Big(\sum_{k=n}^\infty p_k\Big)^{-1}\ \sim\ \dfrac{n^\alpha}{L(n)}
\end{equation}
for some $\alpha\in(0,1)$. Let $T$ be an exponentially distributed
random variable with mean $1$, which is independent of a
subordinator $\{U_t:t\geq 0\}$ with zero drift and L\'{e}vy
measure (\ref{lev}).

It is well known and follows, for example, from our Proposition
\ref{weakN} that $N_n/w(n)$ weakly converges to the Mittag-Leffler
distribution with parameter $\alpha$. From (\ref{momentsQ}) or
from Proposition 3.1 in \cite{Ca} we have
$$
\me \left(\int_0^T e^{-U_t}\,dt\right)^n
\ =\ \dfrac{n!}{\Gamma^n(1-\alpha)\Gamma (1+n\alpha)},
\quad n\in\mn_0,
$$
which means that $\int_0^T e^{-U_t}\,dt$ has Mittag-Leffler
distribution with parameter $\alpha$. Thus,
\begin{equation}\label{univN}
   \dfrac{N_n}{w(n)}\ \dod\ \int_0^T e^{-U_t}\,dt.
\end{equation}
Let $\eta_\alpha$ be a beta-distributed random variable with
parameters $1-\alpha$ and $\alpha$, i.e. with density $x\mapsto
\pi^{-1}\sin(\pi\alpha) x^{-\alpha}(1-x)^{\alpha-1}$,
$x\in (0,1)$. It is well known (see, for example, Theorem 8.6.3 in
\cite{BGT}) that $(1-S_{N_n-1}/n)^\alpha\dod
\eta^\alpha_\alpha$. It can be checked that
$$
\me\eta_\alpha^{n\alpha}
\ =\ \dfrac{\Gamma(\alpha(n-1)+1)}{\Gamma(1-\alpha)\Gamma(\alpha n+1)},
\quad n\in\mn_0.
$$
From (\ref{momM}) it follows that $e^{-U_T}$ has the same moment sequence.
Therefore, since the distribution of $e^{-U_T}$ is concentrated on $[0,1]$,
it coincides with the distribution of $\eta_\alpha^\alpha$. Thus,
\begin{equation}\label{under}
   \left(1-\dfrac{S_{N_n-1}}{n}\right)^\alpha\ \dod\ e^{-U_T}.
\end{equation}
Now we point out a bivariate result generalizing (\ref{univN}) and (\ref{under}).
\begin{assertion}\label{bivar}
   Suppose (\ref{regvar}) holds. Then,
   $$
   w^{-1}(n)(w(n-S_{N_n-1}),N_n)\ \dod\ (e^{-U_T},\int_0^T e^{-U_t}\,dt),
   $$
   where $\{U_t:t\geq 0\}$ is a subordinator with zero drift and L\'{e}vy
   measure (\ref{lev}).
\end{assertion}
\begin{rem}
   Corollary 3.3 in \cite{Port} states that
   \begin{equation}\label{po}
      \left(
         \dfrac{L(n)}{n^\alpha}(N_{n+1}-1),1-\dfrac{S_{N_{n+1}-1}}{n}
      \right)
      \ \dod\ (X,Y),
   \end{equation}
   where the distribution of a random vector $(X,Y)$ was defined by
   the moment sequence. Our proof of Proposition \ref{bivar} is
   different from and simpler than Port's proof of (\ref{po}).
\end{rem}
\begin{proof}
   According to Proposition \ref{mixed} it suffices to verify that
   \begin{equation}\label{inter}
      \lin \dfrac{\me w^i(n-S_{N_n-1})N_n^j}{w^{i+j}(n)}
      \ =\
      \dfrac{j!\,\Gamma(\alpha(i-1)+1)}
      {\Gamma^{j+1}(1-\alpha)\,\Gamma(\alpha(i+j)+1)},
      \quad i,j\in\mn_0.
   \end{equation}
   For $i=0$, Eq. (\ref{inter}) follows from (\ref{libeh}). For $i\in\mn$,
   Eq. (\ref{inter}) is checked as follows.
   \begin{eqnarray*}
      &   & \hspace{-15mm}\me w^i(n-S_{N_n-1})N_n^j\\
      & = & \sum_{k=1}^n\sum_{l=0}^{n-1}
            w^i(n-l)k^j\mmp\{N_n=k,S_{k-1}=l\}\\
      & = & w^i(n)\mmp\{\xi\geq n\}
            + \sum_{l=1}^{n-1}w^i(n-l)\mmp\{\xi\geq n-l\}
            \sum_{k=2}^{l+1}k^j\mmp\{S_{k-1}=l\}\\
      & = & w^i(n)\mmp\{\xi\geq n\} + \sum_{l=1}^{n-1}w^{i-1}(n-l)
            \sum_{k=2}^{l+1}k^j\mmp\{S_{k-1}=l\}.
   \end{eqnarray*}
   As on p.~26 in \cite{Als}, define the function $f(x):=0$ on $[0,1)$
   and $f(x):=(k+1)^j$ on $[k,k+1)$ for $k\in\mn$, and set
   $F(t):=\int_0^t f(x)\,dx$. Then,
   $$
   \sum_{l=1}^{n-1}\sum_{k=2}^{l+1}k^j\mmp\{S_{k-1}=l\}
   \ =\ \sum_{k=1}^{n-1}(k+1)^j\mmp\{N_n>k\}
   \ =\ \me F(N_n).
   $$
   By Karamata's theorem, $F(t)\sim (j+1)^{-1}t^{j+1}$.
   Since $\lim_{n\to\infty}N_n=\infty$ almost surely and
   $(N_n/w(n))^{j+1}\dod\xi_\alpha^{j+1}$, where $\xi_\alpha$ is Mittag-Leffler
   distributed with parameter $\alpha$, we have
   \begin{equation}\label{inter3}
      \dfrac{F(N_n)}{w^{j+1}(n)}\ \dod\ \dfrac{\xi_\alpha^{j+1}}{j+1}.
   \end{equation}
   By (\ref{libeh}), $\lim_{n\to\infty}\me (N_n/w(n))^{j+2}=
   \me\xi_\alpha^{j+2}<\infty$. Therefore, the sequence
   $\{F(N_n)/w^{j+1}(n):n\in\mn\}$ is uniformly integrable which
   together with (\ref{inter3}) implies
   \begin{equation}\label{inter2}
      \me F(N_n)
      \ \sim\ \me \dfrac{\xi_\alpha^{j+1}}{j+1}w^{j+1}(n)
      \ \sim\ \dfrac{j!}{\Gamma^{j+1}(1-\alpha)\Gamma(1+(j+1)\alpha)}
              \dfrac{n^{\alpha(j+1)}}{L^{j+1}(n)}.
   \end{equation}
   Thus, if $i=1$, we have
   $$
   \me w(n-S_{N_n-1})N_n^j
   \ \sim\
   \dfrac{j!}{\Gamma^{j+1}(1-\alpha)\Gamma(1+(j+1)\alpha)}
   \dfrac{n^{\alpha(j+1)}}{L^{j+1}(n)},
   $$
   and (\ref{inter}) follows. Assume now that $i\geq 2$.
   Since $w^{i-1}(n)\sim n^{\alpha(i-1)}/L^{i-1}(n)$,
   Corollary 1.7.3 in \cite{BGT} yields
   $$
   W(s)\ :=\ \sum_{n=1}^\infty s^n w^{i-1}(n)
   \ \sim\ \dfrac{\Gamma(1+\alpha(i-1))}{(1-s)^{1+\alpha(i-1)}L^{i-1}((1-s)^{-1})},
   \quad s\uparrow 1.
   $$
   By the same Corollary, (\ref{inter2}) implies
   \begin{eqnarray*}
      R(s)
      & := & \sum_{n=1}^\infty s^n
            \Big(\sum_{k=2}^{n+1}k^j\mmp\{S_{k-1}=l\}\Big)\\
      & \sim & \dfrac{j!}{\Gamma^{j+1}(1-\alpha)}
            \dfrac{1}{(1-s)^{\alpha(j+1)}L^{j+1}((1-s)^{-1})},
            \quad s\uparrow 1.
   \end{eqnarray*}
   Therefore,
   $$
   W(s)R(s)\ \sim\
   \dfrac{\Gamma(1+\alpha(i-1))j!}{\Gamma^{j+1}(1-\alpha)}
   \dfrac{1}{(1-s)^{1+\alpha(i+j)}L^{i+j}((1-s)^{-1})},
   \quad s\uparrow 1.
   $$
   The sequence $\{w^{i-1}(n):n\in\mn\}$ is non-decreasing. Hence,
   the sequence $\{\sum_{l=1}^{n-1}w^{i-1}(n-l)\sum_{k=2}^{l+1}
   k^j\mmp\{S_{k-1}=l\}:n=2,3,\ldots\}$ is non-decreasing too.
   Another appeal to Corollary 1.7.3 in \cite{BGT} gives, as $n\to\infty$,
   $$
   \sum_{l=1}^{n-1}w^{i-1}(n-l)\sum_{k=2}^{l+1}k^j\mmp\{S_{k-1}=l\}
   \ \sim\
   \dfrac{\Gamma(1+\alpha(i-1))j!}{\Gamma^{j+1}(1-\alpha)\Gamma(1+\alpha(i+j))}
   \dfrac{n^{\alpha(i+j)}}{L^{i+j}(n)}.
   $$
   From this, (\ref{inter}) follows.
\end{proof}
\section{Proof of Theorem \ref{bas1}}
Our proof essentially relies upon the following classical result
\begin{equation}\label{bas2}
   \lin\mmp\{n-S_{N_n-1}\leq k\}
   \ =\ m^{-1}\sum_{i=1}^k\mmp\{\xi\geq i\}\ =:\ \mmp\{W\leq k\},
   \quad k\in\mn.
\end{equation}
In order to see why (\ref{bas2}) holds, note that
\begin{eqnarray*}
   \mmp\{n-S_{N_n-1}=k\}
   & = & \sum_{i=1}^n\mmp\{S_{i-1}=n-k,S_i\geq n\}\\
   & = & \mmp\{\xi\geq k\}\sum_{i=0}^{n-k}\mmp\{S_i=n-k\}\\
   & \to & m^{-1}\mmp\{\xi\geq k\},\quad n\to\infty,
\end{eqnarray*}
by the elementary renewal theorem, and (\ref{bas2}) follows.

From (\ref{impo1}) we conclude that
\begin{equation}\label{impofinite}
   M_n-N_n\ \dod\ M_W'-1,
\end{equation}
where $W$ is a random variable with distribution (\ref{bas2})
which is independent of $\{M_n':n\in\mn\}$. Therefore, for
any sequence $\{d_n:n\in\mn\}$ such that $\lim_{n\to\infty}d_n=\infty$,
\begin{equation}\label{bas4}
   \dfrac{M_n-N_n}{d_n}\ \tp\ 0.
\end{equation}
Assume that the distribution of $\xi$ does not belong to the
domain of attraction of any stable law with index $\alpha\in
[1,2]$. Then, as is well known, it is not possible to find
sequences $x_n>0$ and $y_n\in\mr$ such that $(S_n-y_n)/x_n$
converges to a proper and non-degenerate law. In view of
\begin{equation}\label{equ}
   \mmp\{N_n>m\}\ =\ \mmp\{S_m\leq n-1\},
\end{equation}
the same is true for $N_n$ (see Theorem 7 in \cite{Feller} and/or
Theorem 2 in \cite{Heyde} for more details), and according to
(\ref{bas4}), for $M_n$.

Assume that conditions (ii) of Theorem \ref{bas1} hold.
If $\sigma^2=\infty$ and (\ref{domain}) holds with $\alpha=2$,
then arguing as in the proof of Theorem 2 in \cite{Heyde}
we conclude that, with $a_n$ and $b_n$ defined in our
Theorem \ref{bas1},
$$
\dfrac{N_n-b_n}{a_n}\ \Rightarrow\ \mu_2.
$$
Theorem 5 in \cite{Feller} (if $\sigma^2<\infty$) and
Theorem 7 in \cite{Feller} (if (\ref{domain}) holds for some
$\alpha\in[1,2)$) leads to the same limiting relation (with corresponding
$a_n$ and $b_n$, and with $\mu_2$ replaced by $\mu_\alpha$ in the latter case).

In view of (\ref{bas4}) the same limiting relations hold for $M_n$
and, hence, by Lemma \ref{xm}, for $X_n$. The proof of Theorem
\ref{bas1} is complete.
\section{A probabilistic proof of Theorem \ref{smallerone}}
\label{probab}
Set $Y_n:=n-S_{N_n-1}$. The sequence of distributions of
$\{M_n/\me M_n:n\in\mn\}$ is tight. According to (\ref{moma}),
$\me M_n\sim const\, w(n)$, where $w(n)$ is the same as in
(\ref{regvar2}). Therefore, there exists a sequence $\{n_k:
k\in\mn\}$ such that $\lim_{k\to\infty}n_k=\infty$ and, as $k\to\infty$,
$M_{n_k}/w(n_k)$ converges in law to a random variable $Z$, say,
with a proper law. From $Y_n\tp +\infty$ and the result
of Lemma \ref{repr} we conclude that, as $k\to\infty$,
$\widehat{M}_{Y_{n_k}}/w(n_k)$ converges in law to a
random variable $Z''\od Z$. By Proposition \ref{bivar},
as $k\to \infty$,
$$
\left(
  \dfrac{w(Y_{n_k})}{w(n_k)},\dfrac{N_{n_k}-1}{w(n_k)}
\right)\ \dod\ (M,Q)\ :=\ (e^{-U_T},\int_0^T e^{-U_t}dt).
$$
Rewriting (\ref{impo1}) in the form
$$
\dfrac{M_{n_k}}{w(n_k)}
\ =\
\dfrac{\widehat{M}_{Y_{n_k}}}{w(Y_{n_k})}
\dfrac{w(Y_{n_k})}{w(n_k)}
+\dfrac{N_{n_k}-1}{w(n_k)}
$$
we conclude that, as $k\to\infty$,
$$
\left(
   \dfrac{\widehat{M}_{Y_{n_k}}}{w(Y_{n_k})},
   \dfrac{w(Y_{n_k})}{w(n_k)},
   \dfrac{N_{n_k}-1}{w(n_k)}
\right)\ \dod\ (Z',M,Q),
$$
where $Z'\od Z$ and using characteristic functions it can be
checked that $Z'$ is independent of $(M,Q)$. Furthermore,
\begin{equation}\label{eq2}
   Z\ \od\ MZ'+Q.
\end{equation}
From (\ref{perp}) it follows that the distribution of
$\int_0^\infty e^{-U_t}\,dt$ is a solution of (\ref{eq2}). By
Theorem 1.5 (i) in \cite{Verv} this solution is unique. Therefore,
we have proved that, as $k\to\infty$,
$$
\dfrac{M_{n_k}}{w(n_k)}\ \dod\ \int_0^\infty e^{-U_t}\,dt.
$$
The same argument can be repeated for any sequence like $n_k$,
and the proof is complete.

Combining the proof above with the results of Subsection
\ref{expsubord} immediately give the following corollary.
\begin{cor}
   Suppose (\ref{regvar}) holds. Then,
   \begin{eqnarray*}
      &   & \hspace{-3cm}
            \left(
               \dfrac{M_n-N_n}{w(n-S_{N_n-1})},
               \dfrac{w(n-S_{N_n-1})}{w(n)},
               \dfrac{N_n}{w(n)}
            \right)\\
      & \dod &
      \left(
         \int_0^\infty e^{-(U_{t+T}-U_T)}dt,
         e^{-U_T},
         \int_0^T e^{-U_t}\,dt
      \right).
   \end{eqnarray*}
   Furthermore, $(M_n-N_n)/w(n-S_{N_n-1})$ and
   $(w(n-S_{N_n-1})/w(n),N_n/w(n))$ are asymptotically independent, and
   $$
   w^{-1}_n(M_n-N_n, N_n)\ \dod\
   \left(\int_T^\infty e^{-U_t}\,dt,\int_0^T e^{-U_t}\,dt\right).
   $$
\end{cor}
\section{An analytic proof of Theorem \ref{smallerone}}\label{analytic}
Nothing more than (\ref{basic}) and (\ref{basic2}) is required for
the proof given below. In particular, the construction in Section
\ref{coupl} is not needed.

For $k,n \in \mn$ set $a_k(n):=\me X_n^k$ and $b_k(n):=\me N_n^k$.
For $x\geq 0$ define
$$
\Phi(x)
\ :=\ \frac{\Gamma(1-\alpha)\Gamma(\alpha x+1)}{\Gamma(\alpha(x-1)+1)}-1
\ =\ \alpha x B(\alpha x,1-\alpha)-1,
$$
where $B$ denotes the beta function. Note that
$$
B(\alpha x,1-\alpha)
\ =\ \int_0^1 y^{\alpha x-1}(1-y)^{-\alpha}\,dy
\ =\ \alpha^{-1}\int_0^\infty e^{-xy}(1-e^{-y/\alpha})^{-\alpha}\,dy
$$
and, hence,
\begin{eqnarray}
   \Phi(x)
   & = & \int_0^\infty xe^{-xy}(1-e^{-y/\alpha})^{-\alpha}\,dy-1\nonumber\\
   & = & \int_0^\infty(1-e^{-y/\alpha})^{-\alpha}d(1-e^{-xy})-1\nonumber\\
   & = & \int_0^\infty
         (1-e^{-xy})\frac{e^{-y/\alpha}}{(1-e^{-y/\alpha})^{\alpha+1}}\,dy.
         \label{inte}
\end{eqnarray}
Thus, the function $\Phi$ is the Laplace exponent of an infinitely
divisible law with zero drift and L\'{e}vy measure $\nu$ given in
(\ref{lev}). Note that (\ref{inte}) corrects an error on p.~102 in
\cite{BerYor}. Assuming that (\ref{regvar}) holds we will prove that
\begin{equation}\label{moma}
   \lin \frac{L^k(n)}{n^{\alpha k}}\,a_k(n)
   \ =\ \frac{k!}{\Phi(1)\cdots\Phi(k)}\ =:\ a_k,
   \quad k\in\mn.
\end{equation}
This will imply (see, for example, \cite{BerYor}) that (i) $a_k=
\me(\eta^k)$, $k\in\mn$, where $\eta$ is a random variable with
distribution of the exponential integral of a subordinator with
zero drift and L\'{e}vy measure $\nu$, and that (ii) the moments
$a_1,a_2,\ldots$ uniquely determine the law of $\eta$. Note that
the statement in (i) was first obtained in Example 3.4 in
\cite{Urbanik}. From (i) and (ii) it will follow that (\ref{moma})
implies (\ref{weak}).

Exactly in the same way as for $b_k(n)$ in the proof of Lemma
\ref{DarKacle}, but starting with (\ref{basic}) instead of
(\ref{eqd}), it follows that
$$
a_1(n)\ =\ 1+r_n\sum_{i=1}^{n-1}a_1(n-i)p_i,
$$
and, for $k\in\{2,3,\ldots\}$,
\begin{eqnarray}
   a_k(n)
   & = & D_k(a_1(n),\ldots,a_{k-2}(n))+ka_{k-1}(n)
         +r_n\sum_{i=1}^{n-1}a_k(n-i)p_i\nonumber\\
   & =: & d_k(n)+r_n\sum_{i=1}^{n-1}a_k(n-i)p_i, \label{reca}
\end{eqnarray}
where the $D_k(.)$ are the same as in the proof of Lemma
\ref{DarKacle}, and $r_n:=1/(p_1+\cdots+p_{n-1})$. We are ready to
prove (\ref{moma}). Again, we use induction on $k$. Suppose
(\ref{moma}) holds for $k\in\{1,2,\ldots,m-1\}$. Set
$$
\beta_1\ :=\ \frac{1}{1-b_1}
\quad\mbox{and}\quad
\beta_k:=\frac{1}{b_{k-1}-k^{-1}b_k}\prod_{i=1}^{k-1}\frac{b_{i-1}}{b_{i-1}-i^{-1}b_i},
\quad k\in\{2,3,\ldots\},
$$
where $b_k:=k!/(\Gamma^k(1-\alpha)\Gamma(1+\alpha k))$,
$k\in\mn$, and note that
\begin{equation}\label{useq}
   a_{k-1}-\beta_k(b_{k-1}-k^{-1}b_k)\ =\ 0.
\end{equation}
In the following we exploit an idea given in the proof of
Proposition 3 in \cite{Gne}. Suppose there exists an $\epsilon>0$
such that $a_k(n)>(\beta_k+\epsilon) b_k(n)$ for infinitely many
$n$. It is possible to decrease $\epsilon$ so that the inequality
$a_k(n)>(\beta_k+\epsilon)b_k(n)+c$ holds infinitely often for any
fixed positive $c$. Thus, we can define $n_c:=\inf\{n\geq 1:
a_k(n)>(\beta_k+\epsilon)b_k(n)+c\}$. Then
\begin{equation}\label{gnedin}
   a_k(n)\ \leq\ (\beta_k+\epsilon)b_k(n)+c\mbox{ for all }
   n\in\{1,2,\ldots,n_c-1\}.
\end{equation}
We have
\begin{eqnarray*}
   &   & \hspace{-2cm}(\beta_k+\epsilon)b_k(n_c)+c
   \ < \ a_k(n_c)
   \ \overset{(\ref{reca})}{=}\
         d_k(n_c)+r_{n_c}\sum_{i=1}^{n_c-1}a_k(n_c-i)p_i\\
   & \overset{(\ref{gnedin})}{\leq} &
         d_k(n_c)+c+(\beta_k+\epsilon)r_{n_c}\sum_{i=1}^{n_c-1}b_k(n_c-i)p_i\\
   & \overset{(\ref{reca}),(\ref{recb})}{=} &
   D_k(a)+ka_{k-1}(n_c) + c\\
   &   &  + (\beta_k+\epsilon)(r_{n_c}-1)(b_k(n_c)-D_k(b)-kb_{k-1}(n_c))+\\
   &   &  + (\beta_k+\epsilon)b_k(n_c)-(\beta_k+\epsilon)(D_k(b)+kb_{k-1}(n_c)),
\end{eqnarray*}
or, equivalently,
\begin{eqnarray*}
   0
   & < & D_k(a) + ka_{k-1}(n_c) + (\beta_k+\epsilon)(r_{n_c}-1)(b_k(n_c)-D_k(b)-kb_{k-1}(n_c))\\
   &   & \hspace{1cm} - (\beta_k+\epsilon)(D_k(b)+kb_{k-1}(n_c)),
\end{eqnarray*}
where we have used the abbreviations
$D_k(a):=D_k(a_1(n_c),\ldots,a_{k-2}(n_c))$ and
$D_k(b):=D_k(b_1(n_c),\ldots,b_{k-2}(n_c))$ for convenience.
Divide the latter inequality by
$z(c):=n_c^{(k-1)\alpha}/L^{k-1}(n_c)$ and let $c$ go to $\infty$
(which implies $n_c\to\infty$). Notice that, according to
(\ref{regvar}), $r_n-1\sim n^{-\alpha}L(n)$ and that by the
induction assumption
$$
\lim_{c\to\infty}\frac{D_k(a_1(n_c),\ldots,a_{k-2}(n_c))}{z(c)}\
=\ 0 \quad\mbox{and}\quad
\lim_{c\to\infty}\frac{a_{k-1}(n_c)}{z(c)}\ =\ a_{k-1}.
$$
Using these facts and (\ref{libeh}) we obtain
$$
0\ \leq\
ka_{k-1}+(\beta_k+\epsilon)b_k-(\beta_k+\epsilon)kb_{k-1}.
$$
Since the function $\Phi$ defined at the beginning of the proof
is positive for $x>0$, and $kb_{k-1}/b_k-1=\Phi(k)$, we conclude
that $kb_{k-1}-b_k>0$. Therefore,
$$
\epsilon(kb_{k-1}-b_k)\ \leq\
k(a_{k-1}-\beta_k(b_{k-1}-k^{-1}b_k))\ =\ 0
$$
by (\ref{useq}). This is the desired contradiction. Thus, we have
verified that
$$
\limsup_{n\to\infty}\frac{a_k(n)}{b_k(n)}\ \leq\ \beta_k.
$$
A symmetric argument proves the converse inequality for the lower
bound. Therefore,
$$
a_k(n)\ \sim\ \beta_kb_k(n)\ \sim\ \beta_k b_k
\frac{n^{k\alpha}}{L^k(n)}\ =\ a_k\frac{n^{k\alpha}}{L^k(n)}.
$$
A similar but simpler reasoning yields the result for $k=1$. We
omit the details. The proof is complete.
\section{Proof of Theorem \ref{weaklaw}}
By Lemma \ref{xm} it suffices to prove the result for $M_n$.
Assume first that $m<\infty$. It is well known that
\begin{equation}\label{strong}
   \lin \dfrac{N_n}{n}\ = \ \dfrac{1}{m}
   \quad\text{almost surely.}
\end{equation}
In view of (\ref{impofinite}), $\lim_{n\to\infty}(M_n-N_n)/n=0$ almost
surely, which yields $\lim_{n\to\infty}M_n/n=1/m$ almost surely. By
the elementary renewal theorem, $\me N_n\sim n/m$. Using the same
approach as in Section \ref{analytic} it is straightforward to check that
$\me M_n\sim n/m$.
Conversely, if $M_n/a_n\tp 1$, then (\ref{bas4}) gives
$(M_n-N_n)/a_n\tp 0$. Therefore, $N_n/a_n\tp 1$. An appeal to
(\ref{strong}) allows us to conclude that $a_n\sim n/m$.

Assume now that $m=\infty$. According to (\ref{libeh1}),
$\me N_n^k\sim n^k/L^k(n)$, $k\in\mn$. Again, the same approach
as in Section \ref{analytic} yields
$$
\me X_n^k\ \sim\ \dfrac{n^k}{L^k(n)}\ \sim\ (\me X_n)^k,
\quad k\in\mn.
$$
Therefore,
$$
\lin\me\left(\dfrac{X_n}{\me X_n}\right)^k\ =\ 1,\quad k\in\mn,
$$
which proves (\ref{weaklaw1}). In fact, to arrive at
(\ref{weaklaw1}), it suffices to know that $\me X_n\sim n/L(n)$
and $\me X_n^2\sim n^2/L^2(n)$ and exploit Chebyshev's inequality.
The proof is complete.
\section{Proof of Theorem \ref{equalone}}
By Theorem 3 (c) and formulae on p.~42 in \cite{Bing} (see also
\cite{Haan})
$$
\dfrac{N_n-b(n)-1}{a(n)}\ \Rightarrow\ \mu_1,
$$
where $\mu_1$ is the $1$-stable law with characteristic function
$\int_{-\infty}^\infty e^{itx}\mu_1(dx)=\exp(it\log|t|-|t|\pi/2)$,
$t\in\mr$. By Corollary \ref{joint},
\begin{equation}\label{111}
   \dfrac{M_n}{N_n-1}\ \tp\ 1.
\end{equation}
Therefore,
$$
\dfrac{M_n-b(n)}{a(n)}-\dfrac{M_n-N_n+1}{N_n-1}\dfrac{b(n)}{a(n)}
\ \Rightarrow\ \mu_1.
$$
Thus, to prove the theorem it suffices to show that the
second summand tends to $0$ in probability. Clearly, this can be
regarded as a rate of convergence result for (\ref{111}).
Recalling the notation $Y_n=n-S_{N_n-1}$ and using (\ref{impo1}) gives
\begin{eqnarray*}
   \dfrac{M_n-N_n+1}{N_n-1}\dfrac{b(n)}{a(n)}
   & = & \dfrac{\widehat{M}_{Y_n}}{Y_n/m(Y_n)}\dfrac{m(n)}{m(Y_n)}
         \dfrac{b(n)Y_n}{na(n)}\dfrac{n}{m(n)(N_n-1)}\\
   & =: & \prod_{i=1}^4 K_i(n).
\end{eqnarray*}
By Corollary \ref{joint}, $m(n)M_n/n\tp 1$. Using the equality of
distributions (\ref{impo1}) and the fact that $Y_n\tp\infty$ allows
us to conclude that $K_1(n)\tp 1$. By Theorem 6 in \cite{Er}, $K_2(n)\dod
1/R$, where $R$ is a random variable uniformly distributed on
$[0,1]$. By Proposition \ref{igrek}, $K_3(n)\tp 0$. Finally, by
Corollary \ref{joint}, $K_4(n)\tp 1$. The proof is complete.
\section{Number of collisions in beta coalescents} \label{coal}
In this section the main results presented in Section \ref{intro}
are applied to the number of collisions that take place in
beta coalescent processes until there is just a single block.
Other closely related functionals of coalescent processes such as the
total branch length or the number of segregating sites have been studied
in \cite{DrmIksMoeRoe} and \cite{Moe}.

Let ${\cal E}$ denote the set of all equivalence relations on $\mn$.
For $n\in\mn$ let $\varrho_n:{\cal E}\to{\cal E}_n$
denote the natural restriction to the set ${\cal E}_n$ of
all equivalence relations on $\{1,\ldots,n\}$. For $\eta\in{\cal E}_n$
let $|\eta|$ denote the number of blocks (equivalence classes) of $\eta$.

Pitman \cite{Pit} and Sagitov \cite{Sag} independently introduced
coalescent processes with multiple collisions. These Markovian processes
with state space ${\cal E}$ are characterized by a finite measure
$\Lambda$ on $[0,1]$ and are, hence, also called $\Lambda$-coalescent
processes. For a $\Lambda$-coalescent $\{\Pi_t:t\geq 0\}$, it is known that
the process $\{|\varrho_n\Pi_t|:t\geq 0\}$ has infinitesimal rates
\begin{equation} \label{rates}
   g_{nk}
   \ := \ \lim_{t\downarrow 0} \frac{\mmp\{|\varrho_n\Pi_t|=k\}}{t}
   \ = \ \left(\hspace{-2mm}\begin{array}{c}n\\k-1\end{array}\hspace{-2mm}\right)
         \int_{[0,1]}x^{n-k-1}(1-x)^{k-1}\,\Lambda(dx)
\end{equation}
for all $k,n\in\mn$ with $k<n$. Let $g_n:=\sum_{k=1}^{n-1}g_{nk}$,
$n\in\mn$, denote the total rates. We are interested in the number
$X_n$ of collisions (jumps) that take place in the restricted
coalescent process $\{\varrho_n\Pi_t:t\geq 0\}$ until there is
just a single block. From the structure of the coalescent process
it follows that $(X_n)_{n\in\mn}$ satisfies the recursion
(\ref{basic}), where $I_n$ is independent of $X_2,\ldots,X_{n-1}$
with distribution $\mmp\{I_n=k\}=g_{n,n-k}/g_n$,
$k\in\{1,\ldots,n-1\}$. The random variable $n-I_n$ is the (random)
state of the process $\{|\varrho_n\Pi_t|:t\geq 0\}$ after its
first jump.

We consider beta coalescents, where, by definition, $\Lambda=\beta(a,b)$ is
the beta distribution with density $x\mapsto (B(a,b))^{-1}x^{a-1}(1-x)^{b-1}$ with respect to the Lebesgue
measure on $(0,1)$, and $B(a,b):=\Gamma(a)\Gamma(b)/\Gamma(a+b)$ denotes
the beta function, $a,b>0$. In this case the rates (\ref{rates}) have the form
\begin{eqnarray}
   g_{nk}
   & = & \left(\hspace{-2mm}\begin{array}{c}n\\k-1\end{array}\hspace{-2mm}\right)
         \frac{1}{B(a,b)}
         \int_0^1 x^{a+n-k-2}(1-x)^{b+k-2}\,dx\nonumber\\
   & = & \left(\hspace{-2mm}\begin{array}{c}n\\k-1\end{array}\hspace{-2mm}\right)
         \frac{B(a+n-k-1,b+k-1)}{B(a,b)},
         \quad k,n\in\mn, k<n.
         \label{betarates}
\end{eqnarray}
From
$$
g_{k+1,k}\ =\ \frac{k(k+1)}{2}\frac{B(a,b+k-1)}{B(a,b)}
$$
it follows that
$$
g_n\ =\ \sum_{k=1}^{n-1}(g_{k+1}-g_k)
\ =\ \sum_{k=1}^{n-1}\frac{2}{k+1}g_{k+1,k}
\ =\ \frac{1}{B(a,b)}\sum_{k=1}^{n-1}k B(a,b+k-1).
$$
In the following it is assumed that $b=1$ such that the
rates (\ref{betarates}) reduce to
$$
g_{nk}\ =\ \left(\hspace{-2mm}\begin{array}{c}n\\k-1\end{array}\hspace{-2mm}\right)
\frac{B(a+n-k-1,k)}{B(a,1)}
\ =\ \frac{n!}{(n-k+1)!}a\frac{\Gamma(a+n-k-1)}{\Gamma(a+n-1)},
$$
and the total rates to
$$
g_n\ =\ a\sum_{k=1}^{n-1}kB(a,k)
\ =\
\left\{
   \begin{array}{cl}
      \displaystyle\frac{a}{a-2}
      \Big(1-\frac{\Gamma(a)\Gamma(n+1)}{\Gamma(a+n-1)}\Big)
               & \mbox{for $a>0$, $a\neq 2$,}\\
      2(h_n-1) & \mbox{for $a=2$}.
   \end{array}
\right.
$$
Here, $h_n:=\sum_{i=1}^n 1/i$ denotes the $n$-th harmonic number.
From the last formula it follows that the parameter $a=2$ plays a
special role in this model. Define
\begin{equation} \label{pk1}
   p_k\ :=\
   \frac{(2-a)\Gamma(a+k-1)}{\Gamma(a)\Gamma(k+2)},\quad k\in\mn.
\end{equation}
Assume now that $0<a<2$. In this case (and only in this case) we have
$p_k\geq 0$ for $k\in\mn$ and $\sum_{k=1}^\infty p_k=1$. Let $\xi$ be a
random variable with distribution $\mmp\{\xi=k\}=p_k$, $k\in\mn$.
For $0<a<2$, $a\neq 1$, we can rewrite (\ref{pk1}) in terms of $\alpha:=2-a$
in the form
$$
p_k\ =\ \frac{1}{1-\alpha}{
\left(\hspace{-2mm}\begin{array}{c}\alpha\\k+1\end{array}\hspace{-2mm}\right)
}(-1)^k,\quad k\in\mn.
$$
Therefore, for $a\neq 1$, i.e. $\alpha\neq 1$, $\xi$ has probability generating
function
$$
\me s^\xi
\ =\ \sum_{k=1}^\infty p_ks^k
\ =\ \frac{1}{1-\alpha}\sum_{k=1}^\infty
     \left(\hspace{-2mm}\begin{array}{c}\alpha\\k+1\end{array}\hspace{-2mm}\right)
     (-s)^k
\ =\ \frac{1-\alpha s-(1-s)^\alpha}{(1-\alpha)s}.
$$
For $a=1$, i.e. $\alpha=1$, the probability generating function is
$$
\me s^\xi\ =\ \sum_{k=1}^\infty\frac{s^k}{k(k+1)}
\ =\ 1-\log(1-s)+\frac{\log(1-s)}{s}
$$
with continuous extensions for $s=0$ and $s=1$. For $0<a<2$ it follows by
induction on $n$ that
$$
\mmp\{\xi\geq n\}\ =\
\frac{\Gamma(a+n-1)}{\Gamma(a)\Gamma(n+1)},\quad n\in\mn.
$$
Using $\Gamma(n+x)\sim\Gamma(n)n^x$ for $n\to\infty$, we conclude that
$$
\mmp\{\xi\geq n\}
\ \sim\ \frac{n^{a-2}}{\Gamma(a)}
\ =\ \frac{n^{-\alpha}}{\Gamma(2-\alpha)},\quad n\to\infty.
$$
Thus, if $1<a<2$, or, equivalently, $0<\alpha<1$, Theorem \ref{smallerone}
is applicable (with $L(n)\equiv 1/\Gamma(a)=1/\Gamma(2-\alpha)$),
and we obtain the following result.
\begin{thm} \label{beta1}
   For the $\beta(a,1)$-coalescent with $1<a<2$, i.e., $0<\alpha:=2-a<1$,
   the number $X_n$ of collision events satisfies
   $$
   \frac{X_n}{\Gamma(2-\alpha)n^\alpha}\ \dod \ \int_0^\infty e^{-U_t}\,dt,
   $$
   where $\{U_t:t\geq 0\}$ is a subordinator with zero drift and
   L{\'e}vy measure (\ref{lev}).
\end{thm}
Note that, for $\Lambda=\beta(a,b)$, we have
$\mu_{-1}:=\int x^{-1}\,\Lambda(dx)<\infty$ if and only if $a>1$.
Under the condition $\mu_{-1}<\infty$, limiting results similar to that
presented in the above Theorem \ref{beta1} are known for the
number of segregating sites (see, for example, Proposition 5.1 in \cite{Moe})
for general $\Lambda$-coalescent processes with mutation.

Assume now that $0<a<1$. Then, $m:=\me\xi=1/(1-a)<\infty$. It is straightforward
to verify that
$$
\sum_{k=1}^n k^2p_k\ \sim\ \frac{2-a}{\Gamma(a+1)}n^a,\quad n\to\infty.
$$
In particular, the variance of $\xi$ is infinite. Thus, Theorem
\ref{bas1} is applicable (with $L(n)\equiv (2-a)/\Gamma(a+1)=\alpha/\Gamma(3-\alpha)$,
$C:=1/\Gamma(a)=1/\Gamma(2-\alpha)$,
$b_n:=n(1-a)=n(\alpha-1)$ and $c_n:=n^{1/\alpha}$), and yields the following result.
\begin{thm} \label{beta2}
   For the $\beta(a,1)$-coalescent with $0<a<1$, i.e., $1<\alpha:=2-a<2$,
   the number $X_n$ of collision events satisfies
   $$
   \frac{X_n-n(\alpha-1)}{(\alpha-1)^{(\alpha+1)/\alpha}n^{1/\alpha}}
   \ \Rightarrow\ \mu_\alpha,
   $$
   or, equivalently,
   \begin{equation} \label{beta2asy}
      \frac{X_n-n(\alpha-1)}{(\alpha-1)n^{1/\alpha}}\ \dod\ S_\alpha,
   \end{equation}
   where $\me\exp(itS_\alpha)=\exp(|t|^\alpha(\cos(\pi\alpha/2)+i\sin(\pi\alpha/2){\rm sgn}(t)))$,
   $t\in\mr$.
\end{thm}
Gnedin and Yakubovich \cite[Theorem 9]{GneYak} use analytic methods to
verify the same convergence result (\ref{beta2asy}) for $\Lambda$-coalescents
satisfying $\Lambda([0,x])=Ax^a+O(x^{a+\zeta})$, $x\to 0$, $0<a<1$, $\zeta>
\max\{(2-a)^2/(5-5a+a^2),1-a\}$.

Theorems \ref{beta1} and \ref{beta2} do not cover the
asymptotics of $X_n$ for the Bolt\-hau\-sen-Sznitman coalescent, i.e.
the $\beta(a,b)$-coalescent with $a=b=1$. The limiting behaviour of $X_n$
for the Bolthausen-Sznitman coalescent was studied in \cite{IksMoe},
and follows also from our Theorem \ref{equalone} with
$p_k:=1/(k(k+1))$, $L(n)\equiv 1$, $c(x):=x$, $b(x):=x/\log x+
x\log\log x/(\log x)^2$, and $a(x):=b^2(x)/x\sim x/(\log x)^2$.
Therefore,
the asymptotics of $X_n$ for all $\beta(a,1)$-coalescent processes
with $0<a<2$ is clarified. Unfortunately, our method cannot be
used to treat the asymptotics of $X_n$ for $\beta(a,1)$-coalescent
processes with $a\geq 2$, as in this case the crucial assumption
(\ref{basic2}) is not satisfied.
\section{Possible generalizations} \label{general}
We have studied random recursions (\ref{basic}) under the assumption
that
\begin{equation}\label{basic3}
   I_n\ \dod\ \xi
\end{equation}
with specified rate of convergence (\ref{basic2}). If $\me\xi<\infty$,
this specific rate of convergence (\ref{basic2}) ensures that $X_n$
and $N_n$ have the same limiting behaviour. Under the sole condition
(\ref{basic3}) without any assumption on the speed of convergence
such as (\ref{basic2}), the asymptotics of $X_n$ can differ
significantly from that of $N_n$, even if $\me\xi<\infty$.
Assume for example that $I_2\equiv 1$ and that
$\mmp\{I_n=n-1\}=1-\mmp\{I_n=1\}=1/n$ for $n\geq 3$, or,
equivalently, that $\pi_{n,1}=1-\pi_{n,n-1}=1/n$ for $n\geq 3$. In
this case, (\ref{basic3}) is obviously satisfied with $\xi\equiv
1$. In particular, $S_k\equiv k$, $k\in\mn_0$, and $N_n\equiv n$,
$n\in\mn$. It is straightforward to derive the distribution of
$X_n$. We have $\mmp\{X_n=n-1\}=\prod_{i=2}^n\pi_{i,i-1}=2/n$ and,
for $k\in\{1,\ldots,n-2\}$, $\mmp\{X_n=k\}=\pi_{n-k+1,1}
\prod_{i=n-k+2}^n\pi_{i,i-1}=1/n$. Thus, $X_n/n$ is asymptotically
uniformly distributed on $(0,1)$. In particular, $N_n$ and $X_n$ do
not have a similar limiting behaviour.

It is even more evident that the rate of convergence in (\ref{basic3})
will influence the limiting behaviour of $X_n$, if $\me\xi=\infty$, in
particular, when (\ref{domain}) holds.

For the case $\me\xi=\infty$ we left open the interesting
theoretical problem of finding necessary and sufficient conditions
under which $(X_n-b_n)/a_n$ weakly converges to a proper law.
Theorems \ref{weaklaw}, \ref{smallerone} and \ref{equalone}
are our contribution to the one-sided solution of this problem. To
solve the problem in full generality one should, among others,
understand a weak behaviour of $X_n$ under the assumption
$\sum_{k=n}^\infty p_k\sim 1/L(n)$, where $L$ is some slowly
varying function. It seems that this case is not amenable to the
analysis presented in this work.

We concentrated on $M_n$, the number of jumps
of the process $R^{(n)}:=\{R_k^{(n)}: k\in\mn_0\}$, which is an
interesting generalization of random walks. We think it is of
interest to analyse other functionals of $R^{(n)}$ such as
$M_n^{(i)}:=\#\{k\geq 1: R_k^{(n)}-R_{k-1}^{(n)}=i\}$ for some
fixed $i\in\{0,1,\ldots,n-1\}$, or $T_n:=M_n+M_n^{(0)}=\inf\{k\geq 1:
R_k^{(n)}=n-1\}$.
P. Negadajlov has already
checked that Theorems \ref{weaklaw}, \ref{bas1}, \ref{smallerone} and
\ref{equalone} of the present work remain valid with $X_n\od M_n$
replaced with $T_n$.

\vspace{5mm}

{\bf Acknowledgement.}
The authors thank Alexander Gnedin for fruitful comments and discussions,
in particular, for pointing out an error in Section \ref{general} of the
first version of the manuscript.

\end{document}